\documentclass[draft]{amsart}

\usepackage{amssymb}

\usepackage{url}

\newtheorem{thm}{Theorem}[section]

\newtheorem{prop}[thm]{Proposition}

\theoremstyle{remark}
\newtheorem{rem}[thm]{Remark}

\theoremstyle{definition}

\numberwithin{equation}{section}
\numberwithin{thm}{section}

\begin{document}




\title[Absence of eigenvalues]{Absence of eigenvalues for integro-differential operators with
periodic coefficients}

%

\author{Marius Marinel  Stanescu}
\address{Department of Applied Mathematics\\
 University of Craiova\\
  13 A.I. Cuza Str, Craiova\\
  RO-200396, Romania
}
\email{mamas@mail.md}

\author{Igor Cialenco}
\address{Department of Applied Mathematics\\
Illinois Institute of Technology\\
10 West 32nd Str, Bld E1, Room 208\\
Chicago, IL 60616-3793
}
\email{igor@math.iit.edu}


\thanks{MMS acknowledges partial support from the
Romanian Ministry of Education, Research and Youth (C.N.M.P.), grant PN-CDI-II, No. 3021/2007.
The authors express  their thanks to Prof. P.~Cojuhari for
suggesting the problem and many stimulating conversations, and IC thanks
the participants of Analysis Seminar from the Department of
Mathematics, University of Southern California, for their input
and helpful comments.}


\begin{abstract}
Applying  perturbation theory methods,
the absence of the point spectrum for some nonselfadjoint
integro-differential operators is investigated.  The considered differential operators are
of arbitrary order and act in   either $\mathbf{L}_{p}( \mathbb{R} _{+})$
or $\mathbf{L}_{p}(\mathbb{R} ) \ ( 1\leq p<\infty )$.
As an application of general results, new spectral properties of the perturbed  Hill
operator are derived.
\end{abstract}

\maketitle


\section{Introduction}

The spectral theory of some integro-differential operators
is used to get important theoretical results in theory of neutrons
scattering, plasma oscilations,  quantum physics, mechanics and chaos behavior:
(see for instance classical works by J.~Lehner and G.M.~Wing \cite{LehnerWing56},
E.A.~Catchpole \cite{Catchpole73},  D.~Bohm and E.~Grose \cite{BohmGross49}, N.G.~van Kampen
\cite{Kapmen55}, K.M.~Case \cite{Case78}, and recent survey with application to chaos
behavior by J.~McCaw and B.H.J.~McKeller \cite{McCaw2005}).  From this point of view,
it is important to describe  the operator's spectrum and its components. The
essential spectrum can be easily determined by applying Weyl's type theorems about stability of the
essential spectrum (see T.~Kato \cite{Kato}). However, this is not the case for
other components of the spectrum. For nonselfadjoint operators two fundamental properties are
absence and finiteness of the point spectrum. While these problems look similar, the
methods developed for their study are different.
Some important results on the absence of eigenvalues of differential operators of any order have been
obtained by P.~Cojuhari \cite{Cojuhari93}. Also, P.~Cojuhari and M.M.~Stanescu \cite{CojuhariStanescu03},
\cite{Marius98} studied the same problem for integro-differential operators, with the unperturbed operator being
a differential operator with constant coefficients. The absence of the point spectrum for tridiagonal opeartors
have been investigated by C.G.~Kokologiannaki \cite{Kokologiannaki2003}.

In this paper we will investigate the problem
of absence of the point spectrum for a large class of integro-differential operators.
This operators are generally assumed to be
non-selfadjoint, of any order,  and act in one of the spaces $\mathbf{L}_{p}(\mathbb{R}_{+})$
or $\mathbf{L}_{p}( \mathbb{R} ), \  1\leq p<\infty$. Applying methods from perturbation theory,
we consider the original operator as a sum of a differential operator with periodic coefficients
(the unperturbed operator) and an integro-differential operator (the perturbation).
We establish sufficient conditions on the coefficients and kernels of the perturbation that guarantee
that the point spectrum of the original operator is empty. The paper is organized as follows.
In Section~2 we state the problem and derive some auxiliary results, mainly describing explicitly
the spectrum and the resolvent of the unperturbed operator by applying Floquet theory.
In Section~3 we prove the main result.
The absence of the point spectrum depends on how fast the coefficients and the kernels of the perturbation
decay to zero at infinity. The polynomial decay, with order of decay depending on the multiplicity of
the corresponding Floquet multiplicators, together with subdiagonal property of the kernels
$(k(t,s)=0, \ s<t)$, will guarantee the absence of the eigenvalues of the  perturbed operator.
The results agree with those particular cases established in
\cite{Cojuhari93}, \cite{CojuhariStanescu03}, \cite{Marius98}, and the conditions are
in some sense necessary (see for instance \cite{Cojuhari90}).
In the last section, as an application of general results, we consider the perturbed Hill operator, that represents
and important and interesting results by itself.


\section{The problem and some auxiliary results}
In the space $\mathbf{L}_p(\mathbb{R}_+)$ consider the
differential operator $D=i\frac{d}{dx}$ with the domain of definition determined by the set of all
functions $u\in \mathbf{L}_{p}( \mathbb{R} _{+}) $ which are
absolutely continuous on every bounded interval of the positive semi-axis
and the generalized derivative $u^{\prime }$ belonging to $\mathbf{L}_{p}( \mathbb{R} _{+})$.

Let $H$ be an integro-differential operator of the form

\begin{equation}\label{eq1}
H=\sum\limits_{j=0}^{n}H_{j}D^{j},
\end{equation}
where
%
\[
( H_{j}u) (t) =h_{j}(t) u( t)+\int\limits_{\mathbb{R}
_{+}}k_{j}( t,s) u( s) ds \quad (j=0,\dots,n) \ ,
\]
the functions $h_{j}(t)$ and the kernels
$k_{j}(t,s), \ (j=0,\dots,n; \ t,s\in \mathbb{R}_{+})$,
 are complex-valued functions and smooth as it will be
necessary. We consider the operator $H$ on its maximal domain, i.e.
on the set of all functions $u\in \mathbf{W}_{p}^{n}( \mathbb{R}_{+}) \ (1\leq p<+\infty)$
(where $\mathbf{W}_{p}^{n}( \mathbb{R} _{+})$ denotes the
Sobolev space of order $n$ over $\mathbb{R}_{+}$) such that
$(H_{j}D^{j}) ( u) \in \mathbf{W}_{p}^{n}( \mathbb{R}_{+}) \ (j=0,\dots,n)$.

Assume that the functions $h_j$ have the representation
$h_j(t)= a_j(t) + b_j(t)$ for $t\in\mathbb{R}_+, \ j=0,\dots,n$, such that
$a_j$ are periodic functions of period $T$, $a_j(t+T)=a_j(t)$, and suppose that  $a_n (t)\equiv 1$.
The operator $H$ will be considered as a perturbation of the differential
operator
$A=\sum\limits_{j=0}^{n}A_{j}D^{j}$ by the operator
$B=\sum\limits_{j=0}^{n}B_{j}D^{j}$,
where $A_{j}$ and $B_{j} \ (j=0,\dots,n)$ are operators acting in $\mathbf{L}_p(\mathbb{R}_+)$
and defined by
\[
( A_{j}u) ( t) =a_{j}( t) u( t) \ , \quad
( B_{j}u) ( t) =b_{j}( t) u( t)
+\int\limits_{ \mathbb{R} _{+}}k_{j}( t,s) u( s) ds.
\]

Under above notations, $H=A+B$, where $A$ is a differential operator with periodic coefficients and
$B$ is an integro-differential operator.

{\it The problem is to find sufficient conditions on the coefficients} $b_j$ and kernels
$k_j$, $j=1,\dots,n$, that guarantee that {\it the point spectrum} (the set of all eigenvalues, including
those on the continuous spectrum) of the perturbed operator $H$ {\it is absent}.
To apply perturbation methods from operator theory, we need to
have at hand a manageable representation of the resolvent function $(A-\lambda I)^{-1}$
of the unperturbed operator $A$.

The spectral properties of the operator $A$ have been investigated by many
authors (see for instance  \cite{Mc, Rofe} and the references therein).
In \cite{Rofe} the operator $A$ is
considered in the space $L_2(\mathbb R)$, while in \cite{Mc} in $L_p(\mathbb R)\
(1\leq p\leq \infty)$. In these papers it is shown  that the operator $A$ has a purely
continuous spectrum which coincides with the set of those values $\lambda $
(the zone of relative stability) for which the equation $Au=\lambda u$ has a
non trivial solution, bounded on the whole axis.
Although the spectrum of the operator  $H_0$ is well-known
(see for instance \cite{Mc, Rofe}),  we will present here  a different method for
describing explicitly the resolvent of $A$, which  relies on
 Floquet-Liapunov theory
about linear differential equations with periodic coefficients
(see for instance \cite{Hartman,Y-S}).

Without loss of generality we can assume that $T=1$.

Let us consider the equation
\begin{equation}\label{eqn2}
A\varphi=\lambda\varphi \, ,
\end{equation}
where $\lambda$  is a complex number, or in  vector form
\begin{equation}\label{eqn3}
\frac{dx}{dt}=A(t,\lambda)\,x \, ,
\end{equation}
where
$$
A(t,\lambda)=\left( \begin{array}{cccccc}
0 & 1 & 0 & \dots & 0  & 0 \\
0 & 0 & 1 & \dots & 0  & 0 \\
\dots & \dots & \dots & \dots & \dots   & \dots \\
0 & 0 & 0 & \dots & 0  & 1  \\
\lambda-A_0& -A_1 & -A_2 & \dots & -A_{n-2}& -A_{n-1}
\end{array}\right) , \ \
$$
and $x=(u, Du, \dots, D^{n-1}u)^{t}$.

Denote by $U(t) \ (=U(t,\lambda))$  the matriciant of the equation \eqref{eqn3},
i.e., the matrix which satisfies the following system of differential equations
$$
\frac{dU(t)}{dt}=A(t,\lambda)\,U(t), \ \ U(0)=E_n \, ,
$$
where $E_n$ is $n\times n$ identity  matrix.
The matrix $U(1)$  is called the monodromy matrix of the equation \eqref{eqn3} and
the eigenvalues $\rho_1(\lambda),\dots,\rho_m(\lambda)$ of the matrix $U(1)$ are called
the multiplicators (Floquet multiplicators). Also, we will say that
$U(1)$  is the monodromy matrix and
$\rho_1(\lambda),\dots,\rho_m(\lambda)$   are multiplicators
of the operator $A-\lambda I$.

Consider the matrix $\Gamma =\ln U( 1)$, where $\Gamma$ is
one of the solutions of equation $e^{Y}=U(1)$. Note that $\Gamma$ exists since
the monodromy matrix is nonsingular.  Hence, the matrix $U( t) $
admits the Floquet representation
\begin{equation}\label{eqn4}
U( t) =F( t) e^{t\Gamma },
\end{equation}
where $F(t) $ is a nonsingular, differentiable matrix of period $T=1$.
The change of variables  $x=F(t)y$ in \eqref{eqn3} gives
\begin{equation}\label{eqn5}
\frac{dy}{dt}=\Gamma y,
\end{equation}
where $\Gamma$ depends on $\lambda$ only.
The solution of the Cauchy equation \eqref{eqn3} with initial condition $y(0)=y_{0}$
has the form
\begin{equation}\label{eqn5.1}
y( t) =e^{t\Gamma }y_{0}  \ .
\end{equation}

Let us describe explicitly the structure of matrix $\exp(\Gamma t)$. For this,
we write the matrix $\Gamma $ to its Jordan canonical form,
$\Gamma =SJS^{-1}$, where
$J=\mathrm{diag} [J( 1) ,\dots,J( m)]$,
and $J(\alpha), \ \alpha = 1,\dots,m$, are the
Jordan Canonical blocks
\[
J( \alpha ) = \left(
\begin{array}{ccccc}
\lambda _{\alpha } & 1 & 0 & \dots & 0 \\
0 & \lambda _{\alpha } & 1 & \dots & 0 \\
\dots & \dots & \dots & \dots & \dots \\
0 & 0 & 0 & \dots & 1 \\
0 & 0 & 0 & \dots & \lambda _{\alpha }%
\end{array}
\right) .
\]
Thus
\begin{equation}\label{eqn6}
\exp (Jt)= \mathrm{diag} [ t\exp J( 1), \dots,t\exp J( m)]  \  ,
\end{equation} where
$$
\exp \left( t J( \alpha ) \right) = \exp ( t\lambda _{\alpha })
\left(
\begin{array}{cccc}
1 & t & \dots & \frac{t^{p_{\alpha }-1}}{( p_{\alpha }-1) !} \\
0 & 1 & \dots & \frac{t^{p_{\alpha }-2}}{( p_{\alpha }-2) !} \\
\dots & \dots & \dots & \dots \\
0 & 0 & \dots & 1%
\end{array}
\right) ,
$$
where $p_{\alpha }$ is the dimension of the Jordan block $J( \alpha ), \
\alpha =1,\dots,m$.

From \eqref{eqn5}-\eqref{eqn6}, we conclude that  the components of the general solution $y( t) $ of
\eqref{eqn5} are linear combinations of $\exp ( \lambda _{1}t) , \dots,\exp
( \lambda _{m}t) $ with polynomial coefficients in $t$.

Note that if $\mathrm{Re}(\lambda) >0$, then $\vert t^{k}\exp ( t\lambda
) \vert \to \infty$, for $k=1,2,\dots$, and if $\mathrm{Re}(\lambda) =0$,
then $\vert t^{k}e^{t\lambda }\vert \to
\infty $ for $k=1,2,\dots$ and $\vert t^{k}e^{t\lambda }\vert
\to 1$ for $k=0$.
By spectral image theorem, for each eigenvalue $\lambda _{\alpha }, \ \alpha
=1,\dots,m$,  of the matrix $\Gamma$ the corresponding multiplicator
$\rho _{\alpha }=\exp ( \lambda _{\alpha }), \  \alpha
= 1,\dots,m$, is  in interior, exterior or on the unit circle if
$\mathrm{Re}(\lambda _{\alpha })<0$, $\mathrm{Re}(\lambda _{\alpha })>0$,
or $\mathrm{Re} (\lambda_{\alpha })=0$.

\begin{rem}
The solution $y( t) $ of equation  \eqref{eqn5} belongs to $\mathbf{L}_{p}^n( \mathbb{R} _{+}) $
if the coefficients of the
terms $\exp ( t\lambda _{\alpha })$ with  $\mathrm{Re}\lambda _{\alpha }\geq 0$ are zero.
Thus, if we have  multiplicators inside the unit circle (and only in this
case), then the equations \eqref{eq1} has a nontrivial solutions in
the space $\mathbf{L}_{p}( \mathbb{R} _{+})$, and the inverse operator
$(A-\lambda I)^{-1}$ does not exist.
\end{rem}

Suppose that $\lambda $ is such that all corresponding
multiplicators satisfy the condition $\vert \rho \vert \geq 1$.
Then the inverse operator (possible
unbounded) of  $A-\lambda I$ exists, and to describe its structure, we consider the equation $Au-\lambda
u=\nu $, where $\nu $ is an arbitrary element from $\mathrm{Ran}(A-\lambda I)$.
Similarly to \eqref{eqn3}, we write the last equation in its  vector-form
\begin{equation}\label{eqn7}
\frac{dx}{dt}=A( t,\lambda ) +f,
\end{equation}
where $f=( 0,\dots,\nu ) ^{t}$. The change of variable $x=F( t) y$
in equation \eqref{eqn7} implies
\begin{equation}\label{eqn8}
\frac{dy}{dt}=\Gamma y+F^{-1}( t) f.
\end{equation}
The vector-valued function
$$
y( t) =-\int\limits_{t}^{\infty }\exp ( \Gamma (
t-s) ) F^{-1}( s) f( s) ds
$$
is the solution of nonhomogeneous equation \eqref{eqn8}, and hence
the solution of equation \eqref{eqn7} has the form
\begin{equation}\label{eqn9}
x( t) =-F( t) \int\limits_{t}^{\infty }\exp (
\Gamma ( t-s) ) F^{-1}( s) f( s) ds.
\end{equation}

Taking into account relations \eqref{eqn4}-\eqref{eqn6} and representation
\eqref{eqn9} we get
\begin{equation}\label{eqn10}
\left( ( A-\lambda I) ^{-1}\nu \right) ( t)
=\sum\limits_{\alpha =1}^{m}\sum\limits_{k=0}^{p_{\alpha }}g_{\alpha
k}( t) \int\limits_{t}^{\infty }( t-s) ^{k}\exp (
\lambda _{\alpha }( t-s) ) h_{\alpha k}( s) \nu
( s) ds,
\end{equation}
where $g_{\alpha k}$ and $h_{\alpha k}$ are some continuous and periodic
functions, with period $T=1$.
\begin{rem}\label{rem2}
If the unperturbed operator $A$ acts in $\mathbf{L}_p(\mathbb{R})$, then $\lambda\in\sigma(A)$ if and only if
there exists at least one multiplicator which lie on the unit circle
$\mathbb{T} = \{ z\in\mathbb{C} \ : \   |z| = 1 \}$. Moreover, the point spectrum
of $A$ is absent (for details, see for instance  \cite{IgorPeriodicCoeff2007}).
\end{rem}

\section{The main result}\label{sectionMainResult}

In this section will present some general results about the absence of the point spectrum
of the perturbed operator $H = A + B$. A natural condition, typical for perturbation methods,
is to assume that the perturbation $B$ is subordinated, in some sense, to the unperturbed operator $A$.
In what follows, we assume that $b_n(t) = 0$ and $k_n(t,s) = 0$, for every $t,s\in\mathbb{R}_+$.

By Weyl's type theorem, if the perturbation $B$ is a compact operator, then the essential
spectrum of operators $H$ and $A$ coincide. This is true, for example, if the coefficients
$b_j$'s decay fast enough to zero, as $t\to\infty$, and the kernels $k_j$'s are completely continuous.
However, even if the unperturbed operator $A$ has
no eigenvalues, the operator $H$ can have infinitely many eigenvalues,
including on continuous spectrum. Some more restrictive conditions on the coefficients and kernels will
imply the absence of point spectrum of $H$.

%
%

The following result hold true.

\begin{thm}\label{th1}
 Let $\rho _{\alpha }=\rho
_{\alpha }( \lambda ), \ ( \alpha =1,\dots,m)$
be the Floquet multiplicators corresponding to the operator
$A-\lambda I$ such that $|\rho _{\alpha }| \geq 1
\ ( \alpha =1,\dots,m)$. Assume that $l$ is the maximum
order of canonical Jordan blocks corresponding to
unimodular multiplicators  $|\rho _{\alpha } |=1$.
If there exists $\delta >l$ such that
\[
( 1+t) ^{\delta }b_{j}( t) \in \mathbf{L}_{\infty}( \mathbb{R}_{+}),
\quad j=0,\dots,n \ ,
\]
the integral operators with kernels
\[
( 1+t) ^{\delta }k_{j}( t,s)\, ,  \quad  \delta >l,
j=0,\dots,n\, ,
\]
are bounded in  $\mathbf{L}_{p}( \mathbb{R}_{+})$, and
$$
k_{j}( t,s) =0, \quad (t>s, \ j=0,\dots,n) \ ,
$$
then $\lambda$ is not an eigenvalue of the perturbed operator $H$.
\end{thm}
\begin{proof}
To simplify the presentation of the proof, we will introduce several auxiliary notations.

Denote by $\mathcal{C}$ the Banach space obtained as the
direct sum of $n$ copies of $\mathbf{L}_{p}( \mathbb{R} _{+})$, i.e.
$\mathcal{C}= \oplus_{j=0}^{n-1}\mathbf{L}_{p}( \mathbb{R} _{+})$.
We define the norm in $\mathcal{C}$ as follows
$\Vert \psi \Vert _{\mathcal{C}}=\sum\limits_{j=0}^{n-1}
\Vert \psi _{j}\Vert _{L_{p}( \mathbb{R} _{+}) }$ with
$\psi :=( \psi _{j}) _{j=0}^{n-1}\in \mathcal{C}$.

Let $S$ de the operator acting on $W_{p}^{n}(\mathbb{R} _{+})$ with values in
$\mathcal{C}$, and defined by
\[
Su=( u,Du,\dots,D^{n-1}u) \quad u\in W_{p}^{n}( \mathbb{R}_{+}) \, .
\]
We also consider the following family of operators
\[
( T_{j}u) ( t) =b_{j}( t) u( t)
+\int\limits_{ \mathbb{R} _{+}}k_{j}( t,s) u( s) ds \quad
( t\in \mathbb{R} _{+}, j=0,\dots,n-1)
\]
which, obviously, are bounded in $\mathbf{L}_{p}( \mathbb{R}_{+}) $,
and we associate to this family the operator $T$ acting in the
space $\mathcal{C}$ and defined by
\[
T\psi =\sum\limits_{j=0}^{n-1}T_{j}\psi _{j}\quad ( \psi =( \psi
_{j}) _{j=0}^{n-1}\in \mathcal{C}) \, .
\]
Note that $B=TS$ and $H=A+TS$.

For every  $\tau \geq 0$, we define
$$ 
( L_{\tau }x) ( t) =( 1+t)^{\tau }x(t) \quad t\in\mathbb{R}_+\, ,
$$ 
and for every $p\in[1, \infty)$, we consider the following family of spaces
\[
\mathbf{L}_{p,\tau }( \mathbb{R} _{+}) := \{ u\in \mathbf{L}
_{p}( \mathbb{R} _{+}) \mid L_{\tau }u\in \mathbf{L}_{p}(
\mathbb{R} _{+}) \} .
\]
with corresponding norm $\Vert u\Vert _{p,\tau }:=\Vert L_{\tau }^{-1}u\Vert$.

Suppose by the contrary, that $\lambda $ is an
eigenvalue of $H$, i.e. there exists an element $u\in \mathbf{L}_{p}(
\mathbb{R} _{+}) $, $u\neq 0$, such that
\begin{equation}\label{eqn11}
Hu=\lambda u\, .
\end{equation}
Taking into account that $H=A+TS$, and since $\lambda $ cannot be an
eigenvalue of $A$, the equation \eqref{eqn11} implies
\[
Su+S( A-\lambda I) ^{-1}TSu=0 \, .
\]
We note that $Su\neq 0$, since otherwise the equation \eqref{eqn11} would imply
that $Au=\lambda u$ with $u\neq 0$, that is a contradiction.
In what follows we denote $x=Su$. The equation \eqref{eqn11} written in a vectorial form
gives
\begin{equation}\label{eqn12}
\frac{dx}{dt}=A( t,\lambda ) x+B( t) x,
\end{equation}
where
\[
B( t) = \left(
\begin{array}{cccc}
0 & 0 & \dots & 0 \\
0 & 0 & \dots & . \\
\dots & \dots & \dots & \dots \\
0 & 0 & \dots & . \\
-B_{0} & -B_{1} & \dots & -B_{n-1}
\end{array}
\right) , \ x=\left(
\begin{array}{c}
u \\
Du \\
\dots \\
D^{n-1}u \end{array} \right).
\]
The change of variables $x=F( t) y$ in \eqref{eqn12} implies
\[
\frac{dy}{dt}=\Gamma y+F^{-1}( 1) B( t) y
\]
and consequently, we deduce
\begin{equation}\label{eqn13}
x( t) =-F( t) \int\limits_{t}^{\infty }\exp (
\Gamma ( t-s) ) F^{-1}( s) B( s)
x( s) ds.
\end{equation}
Note that the vector's components from the right hand side of
\eqref{eqn13} are the sums of the following quantities
$$
( K_{\alpha }x_{j}) ( t) =q( t)
\int\limits_{t}^{\infty }( t-s) ^{l_{\alpha }-1}\exp (
\lambda _{\alpha }( t-s) ) p( s) B_{n-j+1}(
s) x_{j}( s) ds \, ,
$$ 
where $p( t) $ and $q( t) $ are continuous periodic
functions of period $1$, $l_{\alpha }$ takes one the values
$1,\dots,p_{\alpha }$, and $\alpha =1,\dots,m, \  j=0,\dots,n-1$.

To complete the proof will use the following result.

Suppose that operators $A$ and $B$ act in a Banach space $\mathcal{D}$, and assume that:
\begin{itemize}
\item[(i)] $\sigma_p(A)=\emptyset$;
\item[(ii)] $B=TS$, with $S$ acting from $\mathcal{D}$ into $\mathcal{C}$, and
$T$ is acting from $\mathcal{C}$ into $\mathcal{D}$, provided that
$\mathrm{Dom} (S)\supset \mathrm{Dom}(T)$;
\item[(iii)] There exists a family of operators $L_{\tau}, \ \tau \geq 0$, on $\mathcal{C}$,
such that for every  $\tau \geq 0$ the operator $L_{\tau }$ is
one-to-one, i.e. $\mathrm{Ker}(L_{\tau }) = 0$. In addition,
$L_{0}=I_{\mathcal{C}}\  (I_{\mathcal{C}}$ is the unit operator in the space $\mathcal{C}$).
\item[(iv)] There exists $\tau \geq 0$ such that if $\psi \in
\mathcal{C}$ and $T\psi \in \mathrm{Ran}( A-\lambda I) $, then
$\psi \in \mathcal{C}_{\tau }$, $S( A-\lambda I) ^{-1}T\psi \in
\mathcal{C}_{\tau }$ and
\[
\Vert S( A-\lambda I) ^{-1}T\psi \Vert _{\mathcal{C}
,\tau }\leq a\Vert \psi \Vert _{\mathcal{C},\tau } \quad (0<a<1) \ ,
\]
where $|u|_\tau := \| L_\tau u\|_|{\mathcal{C}}$ for $u\in D_\tau:= \mathrm{Dom}(L_\tau)$;
\item[(v)] For every $\psi \in \mathcal{C}_{\tau }$ such that
$T\psi \in \mathrm{Ran}( A-\lambda I) $, the following inequality
holds true
\[
\Vert S( A-\lambda I) ^{-1}T\psi \Vert _{\mathcal{C}
,\tau }\leq c\Vert \psi \Vert _{\mathcal{C},\tau ^{\prime}} \, ,
\]
where $\tau >\tau ^{\prime }\geq 0$  and $c$ is a positive constant independent of $\psi$.
\end{itemize}
\vspace{-.1cm}
Conditions (i)-(v) imply that $\lambda$ is not an eigenvalue of the perturbed operator $A+B$.
For detailed proof see for instance \cite{Cojuhari93}.

Following the same notations, we observe that our operators satisfy conditions (i)-(iii).
To check the conditions (iv) and (v) we define the operator
\[
( R( \lambda ) u) ( t)
=\int\limits_{t}^{\infty }\exp ( \lambda ( t-s) )
u( s) ds \quad  ( 0<t<\infty ) \, .
\]
For all $\tau \geq 0$ and $\mathrm{Re}(\lambda) >0$ the operator
$L_{\tau}^{-1}R( \lambda ) L_{\tau }$ is bounded in $\mathbf{L}_{p}(
\mathbb{R} _{+})$, since $\Vert L_{\tau }^{-1}R( \lambda
) L_{\tau }\Vert \leq ( \mathrm{Re}(\lambda) ) ^{-1}$
(see for instance Lemma 1 and 2 from \cite{Cojuhari92}).
Moreover, for all $\varepsilon >0$ we have $a( \tau ) =\Vert
L_{\tau }^{-1}R( \lambda ) L_{\tau +\varepsilon }\Vert
\to 0$, when $\tau \to \infty$. If $\mathrm{Re}\lambda =0$, the
operator $L_{\tau }^{-1}R( \lambda ) L_{\tau +1}$ is bounded in
$\mathbf{L} _{p}( \mathbb{R} _{+})$, given that
$\Vert L_{\tau}^{-1}R( \lambda ) L_{\tau +1}\Vert \leq 2$.
Note that
\[
( R^{m}( \lambda ) u) ( t) =(
-1) ^{m-1}\int\limits_{t}^{\infty }( t-s) ^{m-1}\exp (
\lambda ( t-s) ) u( s) ds\, .
\]

Let us estimate the norm
$$
\Vert ( K_{\alpha }x_{k}) ( t) \Vert
_{p,\tau }=\Vert L_{\tau }^{-1}( K_{\alpha }x_{k}) (
t) \Vert \, .
$$
For $\mathrm{Re}(\lambda _{\alpha })>0$, using the assumptions on
the functions $b_{j}$ and kernels $k_{j} \ (j=0,\dots,n)$,
we obtain the following estimate
\begin{equation}\label{eqn14}
\Vert ( K_{\alpha }x_{k}) ( t) \Vert
_{p,\tau }=ca( \tau ) \Vert x_{k}\Vert _{p,\tau }.
\end{equation}
For $\mathrm{Re}(\lambda _{\alpha }) =0$ the following equalities hold true
$$
\begin{array}{c}
( I-iD) R( 1) x=x, \\
( I-iD) R( \lambda ) x=x+( 1-\lambda )
R( \lambda ) x\, ,

\end{array}
$$
where  $x\in \mathrm{Dom}( R( \lambda )$. The above, together with initial assumptions, implies
\eqref{eqn14}. Hence, \eqref{eqn14} is satisfied for every $\lambda_\alpha$. Consequently, we get
\begin{equation}\label{eqn15}
\Vert x\Vert _{\mathcal{C},\tau }\leq ca( \tau )
\Vert x\Vert _{\mathcal{C},\tau }.
\end{equation}
Similarly to \eqref{eqn15}, for $\mathrm{Re}(\lambda _{\alpha })\geq 0$ we
obtain
\begin{equation}\label{eqn16}
\Vert ( K_{\alpha }x_{k}) ( t) \Vert
_{p,\tau }=c( \tau ) \Vert x_{k}\Vert _{p,\tau
^{\prime }} \ ,
\end{equation}
and thus
\begin{equation}\label{eqn17}
\Vert x\Vert _{\mathcal{C},\tau }\leq c( \tau )
\Vert x\Vert _{\mathcal{C},\tau ^{\prime }}
\end{equation}
where $c$ is a  constant, and $\tau >\tau ^{\prime }\geq 0$.

Let $\tau ^{\prime }=\tau -\varepsilon, \ \varepsilon >0$. From
estimate \eqref{eqn17}, it follows that $\Vert x\Vert _{\mathcal{C},
\varepsilon }<\infty$. Hence, again from \eqref{eqn17}, we get
$\Vert x\Vert _{\mathcal{C},2\varepsilon }<\infty $ and, in general
$\Vert x\Vert _{\mathcal{C},n\varepsilon }<\infty $.
Since $\varepsilon$ can be chosen arbitrarily, we have that
$\Vert x\Vert_{ \mathcal{C},\varepsilon }<\infty$
for every $\tau \geq 0$. However, as we mentioned above
$x\neq 0$, and by \eqref{eqn16}, we get $1\leq c( \tau )$.
This is a contradiction, since  $c( \tau ) \to 0$, as $\tau \to \infty $.
This completes the proof.
\end{proof}

We conclude this section with the case of whole real line. Using Remark \ref{rem2},
by similar arguments as in Theorem \ref{th1}, one can prove the following
\begin{thm}\label{th2}
Assume the operator $H$ acts in the
space $L_{p}( \mathbb{R})$,  and $\rho _{\alpha }=\rho _{\alpha }( \lambda )$ are all unimodular
multiplicators $\alpha =1,\dots,m$. Suppose that $l$ is  the maximum value for the orders of
canonical Jordan blocks corresponding to the multiplicators
$\rho _{\alpha} \ ( \alpha =1,\dots,m)$. If
$( 1+\vert t\vert ) ^{\delta }b_{j}( t) \in
\mathbf{L}_{\infty ,\delta }( \mathbb{R}_{+}) \quad ( \delta >l, j=0,\dots,n)$
and the kernels $k_{j}( t,s) \ (j=0,\dots,n)$ are such that $k_{j}(t,s)=0$
for $\vert t\vert >\vert s\vert$, and the integral operators with kernels
$( 1+\vert t\vert ) ^{\delta }k_{j}( t,s) \quad ( \delta >l, j=0,\dots,n)$
are bounded in the space $\mathbf{L}_{p}( \mathbb{R})$,
then $\lambda $ is not an eigenvalue of the operator $H$.
\end{thm}

\section{Application}\label{sectionApplication}

In this section we will apply the general results from Section \ref{sectionMainResult} to perturbed Hill
operator.

In the space $\mathbf{L}_{p}( \mathbb{R} _{+}) \ ( 1\leq p<\infty )$
we consider the following integro-differential operator
\begin{eqnarray*}
( Hu) ( t)  & = & ( D^{2}u) ( t) +p( t) u( t) +b_{1}( t) ( Du)
                ( t) +b_{2}( t) u( t) + \\
            & + &\int\limits_{0}^{\infty }k_{1}( t,s) ( Du)
            (s) ds+\int\limits_{0}^{\infty }k_{2}( t,s) u( s)ds \, ,  \\
            & \ &0 <t <\infty, u\in W_{p}^{2}( \mathbb{R} _{+}),
\end{eqnarray*}
where $p( t+1) =p( t) $, $b_{j}( t) \in
\mathbf{L}_{\infty }( \mathbb{R} _{+}), \ j=1,2$,
and kernels $k_{j}( t,s) $ $\in \mathbf{L}_{\infty }( \mathbb{R}
_{+}\times \mathbb{R} _{+}), \ j=1,2$.

The unperturbed operator
\[
( Au) ( t) =( D^{2}u) ( t)
+p( t) u( t)
\]
is Hill operator (see for example \cite{Glazman}).
It is known (see for instance \cite{Rofe} or \cite{Glazman}) that the multiplicators
corresponding to $\lambda\in\sigma(A)$ are
simple and of modulus $1$. Hence, by Theorem \ref{th1}, we have the following result.

\begin{prop}
If
\[
( 1+t) ^{\delta }b_{j}( t) \in \mathbf{L}_{\infty
,\delta }( \mathbb{R} _{+}) ( \delta >1, j=1,2) ,
\]
the kernels $k_{j}( t,s) \ (j=1,2) $ are such that $k_{j}( t,s) =0$
for $t>s$ and the integral operators with kernels
\[
( 1+t) ^{\delta }k_{j}( t,s) ( \delta >1; j=1,2;
t,s\in \mathbb{R} _{+})
\]
are bounded on the space $L_{p}( \mathbb{R} _{+})$,
then the inner point of the continuous spectrum of the operator
$H $ is not eigenvalue.
\end{prop}

If $\lambda $ is one extreme point of the $H$ operator's continuous spectrum
, then there is only one multiplicator equal to $1$ or $-1$ (see \cite{Rofe} or
\cite{Glazman}). This multiplicator is two-fold. Thus, based on Theorem 1 we obtain
the following statement.

\begin{prop}
If
\[
( 1+t) ^{\delta }b_{j}( t) \in \mathbf{L}_{\infty}( \mathbb{R} _{+}) \quad  \delta >2, j=1,2\,  ,
\]
the kernels $k_{j}(t,s), \ j=1,2$,
are such that $k_{j}( t,s) =0, \ t>s$, and the integral operators with kernels
\[
( 1+t) ^{\delta }k_{j}( t,s) \quad  \delta >2; j=1,2;\ t,s\in \mathbb{R} _{+}\, ,
\]
are bounded on the space $L_{p}(\mathbb{R}_{+})$,
then the extreme points of the $H$ operator's continuous spectrum  cannot be eigenvalues.
\end{prop}

Similar results can be proved for the case when operator $H$ is considered along the whole axis
$\mathbb{R}$. Some particular cases of integro-differential operators (second and forth order) are discussed in
\cite{Almamedov68} and \cite{Almamedov88}.

\bibliographystyle{amsplain}

\begin{thebibliography}{99}


\bibitem{Almamedov68}
M.~S. Almamedov, \emph{The spectrum of a linear integro-differential operator
  on the entire axis}, Izv. Akad. Nauk Azerba\u\i d\v zan. SSR Ser. Fiz.-Tehn.
  Mat. Nauk \textbf{1968} (1968), no.~3, 110--117.

\bibitem{Almamedov88}
\bysame, \emph{Spectrum of a fourth-order linear integro-differential
  operator}, Dokl. Akad. Nauk SSSR \textbf{299} (1988), no.~3, 525--529.

\bibitem{BohmGross49}
D.~Bohm and E.P. Gross, \emph{Theory of plasma oscillations, a origin of medium
  like behavior}, Phys. Rev. \textbf{75} (1949), no.~3, 1185.

\bibitem{Case78}
K.~M. Case, \emph{Plasma oscillations}, Phys. Fluids \textbf{21} (1978), no.~2,
  249--257.

\bibitem{Catchpole73}
E.~A. Catchpole, \emph{An integro-differential operator}, J. London Math. Soc.
  (2) \textbf{6} (1973), 513--523.

\bibitem{IgorPeriodicCoeff2007}
Ig. Cialenco, \emph{On the spectrum of the perturbed differential operators
  with periodic coefficients}, submited for publication,
  \url{http://arxiv.org/abs/0708.0854}.

\bibitem{Cojuhari90}
P.~Cojuhari, \emph{The absence of eigenvalues in a perturbed discrete
  {W}iener-{H}opf operator}, Izv. Akad. Nauk Moldav. SSR Mat. (1990), no.~3,
  26--35, 78.

\bibitem{Cojuhari92}
\bysame, \emph{On the point spectrum of the perturbed {W}iener-{H}opf integral
  operator}, Mat. Zametki \textbf{51} (1992), no.~1, 102--113.

\bibitem{Cojuhari93}
\bysame, \emph{On the spectrum of singular nonselfadjoint differential
  operators}, Operator extensions, interpolation of functions and related
  topics (Timi\c soara, 1992), Oper. Theory Adv. Appl., vol.~61, Birkh\"auser,
  1993, pp.~47--64.

\bibitem{CojuhariStanescu03}
P.~Cojuhari and M.~Stanescu, \emph{Absence of eigenvalues for
  integro-differential operators}, Opuscula Math. (2003), no.~23, 5--14.

\bibitem{Glazman}
I.~M. Glazman, \emph{Direct methods of qualitative spectral analysis of
  singular differential operators}, Translated from the Russian by the IPST
  staff, Israel Program for Scientific Translations, Jerusalem, 1965, 1966.

\bibitem{Hartman}
P.~Hartman, \emph{Ordinary differential equations}, Classics in Applied
  Mathematics, vol.~38, Society for Industrial and Applied Mathematics (SIAM),
  Philadelphia, PA, 2002.

\bibitem{Kato}
T.~Kato, \emph{Perturbation theory for linear operators}, Die Grundlehren der
  mathematischen Wissenschaften, Band 132, Springer-Verlag New York, Inc., New
  York, 1966.

\bibitem{Kokologiannaki2003}
C.~G. Kokologiannaki, \emph{Absence of the point spectrum in a class of
  tridiagonal operators}, Appl. Math. Comput. \textbf{136} (2003), no.~1,
  131--138. \MR{MR1935603 (2003i:47029)}

\bibitem{LehnerWing56}
J.~Lehner and G.~M. Wing, \emph{Solution of the linearized {B}oltzmann
  transport equation for the slab geometry}, Duke Math. J. \textbf{23} (1956),
  125--142.

\bibitem{McCaw2005}
J.~McCaw and B.~H.~J. McKellar, \emph{Pure point spectrum for the time
  evolution of a periodically rank-{$N$} kicked {H}amiltonian}, J. Math. Phys.
  \textbf{46} (2005), no.~3, 032108, 24. \MR{MR2125557 (2005k:81071)}

\bibitem{Mc}
D.~McGarvey, \emph{Operators commuting with translation by one. {I}.
  {R}epresentation theorems}, J. Math. Anal. Appl. \textbf{4} (1962), 366--410.

\bibitem{Rofe}
F.~S. Rofe-Beketov, \emph{On the spectrum of non-selfadjoint differential
  operators with periodic coefficients}, Dokl. Akad. Nauk SSSR \textbf{152}
  (1963), 1312--1315.

\bibitem{Marius98}
M.~M. St{\u{a}}nescu, \emph{On the spectrum of some integro-differential
  operators}, Bul. Acad. \c Stiin\c te Repub. Mold. Mat. (1998), no.~3, 21--28.

\bibitem{Kapmen55}
N.G. van Kapmen, \emph{The theory of stationary waves in a plasma}, Phys
  \textbf{21} (1955), 949.

\bibitem{Y-S}
V.~A. Yakubovich and V.~M. Starzhinski{\u\i}, \emph{Linear differential
  equations with periodic coefficients and their applications}, Izdat.
  ``Nauka'', Moscow, 1972.

\end{thebibliography}
%
%


\def\cprime{$'$} \def\cprime{$'$} \def\cprime{$'$} \def\cprime{$'$}
  \def\cprime{$'$} \def\cprime{$'$} \def\cprime{$'$} \def\cprime{$'$}
  \def\cprime{$'$} \def\cprime{$'$} \def\cprime{$'$} \def\cprime{$'$}
  \def\cprime{$'$} \def\cprime{$'$} \def\cprime{$'$} \def\cprime{$'$}
  \def\cprime{$'$} \def\cprime{$'$}
\providecommand{\bysame}{\leavevmode\hbox to3em{\hrulefill}\thinspace}
\providecommand{\MR}{\relax\ifhmode\unskip\space\fi MR }
\providecommand{\MRhref}[2]{%
  \href{http://www.ams.org/mathscinet-getitem?mr=#1}{#2}
}
\providecommand{\href}[2]{#2}

\end{document}